# Market-based Control of Air-Conditioning Loads with Switching Constraints for Providing Ancillary Services

Yao Yao, Yizhi Cheng, Peichao Zhang
Department of Electrical Engineering
Shanghai Jiao Tong University
Shanghai, China
yaoyatututu@foxmail.com

*Abstract*—Air-conditioning loads (ACLs) are among the most promising demand side resources for their thermal storage capacity and fast response potential. This paper adopts the principle of market-based control (MBC) for the ACLs to participate in the ancillary services. The MBC method is suitable for the control of distributed ACLs because it can satisfy diversified requirements, reduce the communication bandwidth and protect users' privacy. The modified bidding and clearing strategies proposed in this paper makes it possible to adjust the switching frequency and strictly satisfy the lockout time constraint for mechanical wear reduction and device protection, without increasing the communication traffic and computational cost of the control center. The performance of the ACL cluster in two typical ancillary services is studied to demonstrate the effect of the proposed method. The case studies also investigate how the control parameters affect the response performance, comfort level and switching frequency.

*Index Terms*—air-conditioning loads, market-based control, lockout time, switching frequency, ancillary service

## I. INTRODUCTION

The intermittence and volatility of renewable energy resources have posed significant threats on the stability and security of power systems. More regulation resources are required to provide ancillary services to ensure the stable operation.

The air-conditioning load (ACL) has become one of the most promising demand resources to provide ancillary services due to its thermal storage capability, large energy capacity and quick response potential. Many researchers have proposed effective control methods of ACLs in this aspect [1-6].

Considering the large scale of ACLs, the following control problems are identified:

(1) The diversified requirements of the users and the ACLs, e.g. the temperature range and the switching constraints, should be considered and satisfied with an efficient method.

(2) The amount of information that ACLs send to the control center should be minimized, while remaining sufficient, to relieve the communication and computation burden.

(3) Considering the implementation and communication cost, the downlink control should be designed to avoid sending individual control signal to each ACL.

Unfortunately, the existing method can only address some of the above problems. In [1-3], although the comfort requirements can be satisfied, the operation state of each ACL has to be specified and the lockout time cannot be guaranteed. In [4-5], considering lockout time constraints in the model predictive control (MPC) scheme complicates the original problem, and more information is required by the control center in real time. The market-based control (MBC) method in our previous work [6] requires little information from the distributed ACLs and simplifies the downlink control, but has not considered the switching frequency and the lockout time constraint.

To solve all the identified problems, this paper adopts the MBC control framework in [6] because it can satisfy different requirements, has relatively low communication bandwidth and can provide privacy protection. The control strategies are modified to include the additional switching constraints in an efficient way.

Compared with our previous work [6], the contributions of this paper are threefold.

1) New bidding and market-clearing strategies are proposed to consider the factor of switching frequency and lockout time constraint without increasing the implementation cost and computational cost.

2) The parametric studies are conducted to investigate the impact of the control parameters on the response performance, comfort level and switching cycles.

3) A case study based on the frequency regulation ancillary service is conducted to demonstrate the performance of the ACL cluster in the fast response application scenario.

The basic control framework in this paper mainly consists of two parts. The first part is to calculate the target signal of the ACL cluster in different ancillary services. The second part is to allocate the target power to individual ACLs with the MBC method. The constraints are included in the bidding and clearing strategy.

This paper is organized as follows: Section II introduces the typical ancillary services of ACLs and the calculation of the corresponding target power. Section III introduces the allocation strategy of the aggregated power. Section IV shows the simulation results. The conclusion and future work are summarized in Section V.

## II. TYPICAL ANCILLARY SERVICES OF ACLS

Similar to the traditional energy storage devices, ACLs, which have the thermal storage capacity, can shift the power consumption and store the thermal energy in the ambient environment. Considering the comfort requirements, the air temperature can't exceed the allowed temperature range, for which they are resources with limited capacity. To lower the required capacity, it's better for the ACLs to respond to the control signal with average value near zero. The typical ancillary services involved in this paper are the mitigation of microgrid tie-lie power fluctuations and the frequency regulation.

### A. Mitigation of Microgrid Tie-line Power Fluctuations

Microgrid is an aggregation of distributed micro-sources, energy storage and local loads. It involves a large number of renewable resources such as wind power and photovoltaic power generators. To integrate the microgrids into the grid reliably and stably, it's important to mitigate the fluctuations resulted from the intermittent resources. Because of the energy-limited characteristics and the quick-responding ability, the ACLs are suitable to mitigate the high frequency components of the fluctuations.

The microgrid system in [6] is adopted in this paper. Define $P_w$ as the wind power, $P_{AC}$ as the aggregated power of ACLs, $P_L$ as the total power of the uncontrollable loads, and $P_g$ as the tie-line power. Ignoring the line loss, the power balance equation at time $k$ can be written as follows:

$$P_g[k] = P_{AC}[k] + P_L[k] - P_W[k]. \quad (1)$$

The power consumption of the ACLs without external control is called the baseline load, which is denoted by $P_{ACbase}$. According to formula (1), the tie-line power with all ACLs uncontrolled can be calculated as:

$$P_{g0}[k] = P_{ACbase}[k] + P_L[k] - P_W[k]. \quad (2)$$

Based on the discrete low-pass filter principle, the smoothed tie-line power can be calculated in the recursive form:

$$P_{gLPF}[k] = \alpha P_{gLPF}[k-1] + (1-\alpha) P_{g0}[k] \quad (3)$$

where $\alpha = \tau/(\tau + \Delta t)$ is the filter coefficient, $\tau$ is the time constant, and $\Delta t$ is control cycle.

To track the smoothed tie-line power, the target power of the ACLs for the fluctuation mitigation service can be calculated as:

$$P^*_{AC,FM}[k] = P_{ACbase}[k] + P_{gLPF}[k] - P_{g0}[k]. \quad (4)$$

### B. Frequency Regulation

To track the scheduled power and reduce the system frequency deviation, the control center calculates the area control error (ACE) signal and sends it to the ACE resources. PJM decomposes the ACE signal into low-frequency component regA and high-frequency component regD through a low-pass filter. Fast response resources can respond to regD signal to achieve higher payments. The average value of regD signal is zero with in a certain time interval, so the ACLs have great potential to participate in the scheme.

The target power of ACLs for the frequency regulation service can be calculated as:

$$P^*_{AC,REG}[k] = P_{ACbase}[k] - P_{reg}[k] \quad (5)$$

where $P_{reg}$ is the regulation signal which equals to the raw normalized regD signal multiplied by the hourly contracted regulation capacity of the ACLs. The operator is "minus" because in this paper, the generator criterion is adopted, which means when the regulation signal is positive, the ACL cluster has to reduce the comsumed power.

The baseline load is necessary to the control of ACLs. Many papers have studied the estimation method of the baseline load[7-8]. As it is not the key point in this paper, the simplified estimation method used in [6] is adopted.

## III. ALLOCATION METHOD CONSIDERING SWITCHING FREQUENCY

This section introduces the control method to allocate the target power to individual ACLs, and takes the switching frequency and lockout time constraints into consideration.

The MBC framework has three major advantages: 1) The control center does not need to collect the users' comfort requirements and the thermal parameters, which can reduce the amount of the exchanged information and protect the privacy of the users. 2) The control center only needs to broadcast the clearing result as the control signal. It significantly simplifies the downlink control and is suitable for the control of distributed resources. 3) It features low computational cost. And the bid strategy and the clearing strategy can be modified to include additional constraints without increasing the complexity of the problem.

### A. General Control Framework

The general control framework comprises three main stages in each control cycle:

1) Bid Stage: Each ACL has a local controller as its agent. The controller determines the bid information based on local requirements, e.g. the comfort preferences and the switching constraints, and sends the bid information to the control center before the next control cycle.

2) Aggregation Stage: First, the control center calculates the target power $P^*_{AC}$. Then, the virtual market set in the control center aggregates the bids from the ACLs and forms the demand curve. The supply curve is determined based on the target power. Finally, the virtual market solves the intersection of the demand and supply curve and finishes market clearing.

3) Disaggregation stage: The ACLs receive and respond to the broadcast clearing signal $p^*$. In this stage, the target power of the ACL cluster is allocated to individual ACLs.

The specific methods will be introduced in the following three subsections.

### B. Bidding Strategy of the ACLs

The bid information of the $i^{th}$ ACL at the $k^{th}$ control cycle is:

$$B_i[k] = ([p_{bid}, q_{bid}], s)_i[k] \quad (6)$$

where the bid price $p_{bid}$ is only a control signal without economic meaning, bid quantity $q_{bid}$ is the power when the

ACL is on which is normally set as the rated power[6] and $s$ denotes the operation state which equals 1 when the ACL is on and 0 when it's off.

The bid price $p_{bid}$ reflects how much the ACL wants to be switched on in the next control cycle. If the bid price is high, it's more possible for the ACL to turn on; if the bid price is low, it's more possible to turn off.

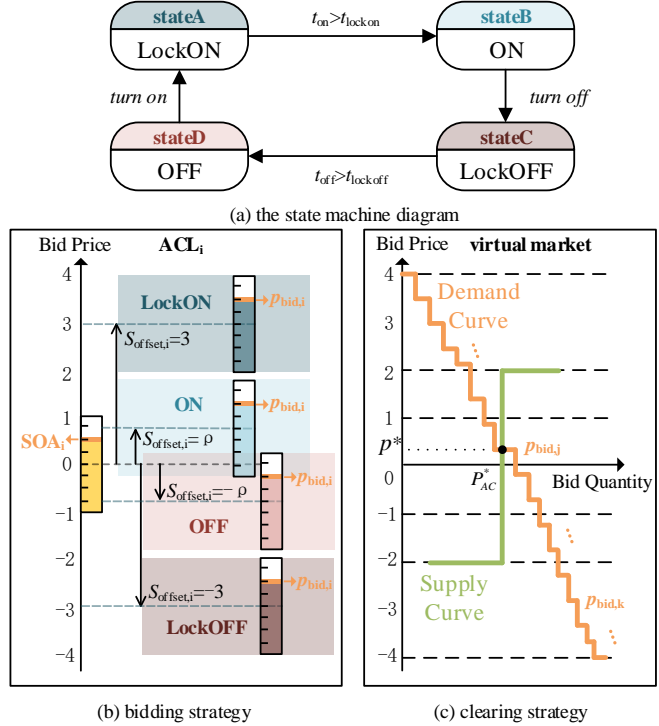

Fig. 1 The bidding strategy and the clearing strategy

Fig. 1(a) shows the state machine diagram of the ACLs, where $t_{on}/t_{off}$ is the duration time since last time the ACL turns on/off, $t_{lockon}/t_{lockoff}$ is the on/off lockout time. When the ACL is in the "LockON/LockOFF" state, it must keep on/off. When the lockout time constraint has been satisfied, the ACL turns into the "ON/OFF" state and is allowed to be switched.

Fig. 1(b) illustrates how the bid price of each ACL is calculated. In this paper, the bid price $p_{bid}$ consists of the following two terms:

1) the state of air-temperature (SOA) which reflects the comfort level of the users and can be calculated as:

$$SOA = \begin{cases} \dfrac{T_{air}-T_{desired}}{T_{max}-T_{desired}}, T_{air} \geq T_{desired} \\ \dfrac{T_{air}-T_{desired}}{T_{desired}-T_{min}}, T_{air} < T_{desired} \end{cases} \quad (7)$$

where $T_{desired}$ is the desired indoor air-temperature, $T_{max}$ and $T_{min}$ is the upper and lower temperature limit when participating in the ancillary services, and $T_{air}$ is the current indoor air-temperature.

It's obvious that $SOA \in [-1,1]$. When it equals to zero, the user's comfort requirement is best satisfied. The higher SOA is, the higher ACL's urgency to be switched on is. The lower SOA is, the higher ACL's urgency to be switched off is.

2) the offset from SOA ($S_{offset}$) which reflects the practical considerations

The first consideration is the switching frequency. Involving the ACLs into ancillary services may increase the ON/OFF cycle times. The bid strategy should be designed to avoid high switching frequency.

The second consideration is the lockout time constraint. It's a rigid constraint according to which the minimal ON/OFF time of each ACL must be guaranteed.

Combined with the above two terms, the bid price can be calculated as:

$$p_{bid,i} = SOA_i + S_{offset,i} \quad (8)$$

$$S_{offset,i} = \begin{cases} 3, stateA : LockON \\ \rho, stateB : ON \\ -\rho, stateD : OFF \\ -3, stateC : LockOFF \end{cases} \quad (9)$$

where $\rho$ is a variable which varies in [0,1].

It's obvious that $S_{offset}$ is determined by the state of the ACL. The mechanism is explained as follows:

When the ACL is in "ON" state or "OFF" state, the operation state of the ACL can be switched freely. The switching frequency can be qualitatively adjusted by the variable $\rho$. If $\rho=1$, $p_{bid}$ of the "ON" ACL varies in [0,2] and $p_{bid}$ of the "OFF" ACL varies in [-2,0]. As a result, all the "ON" ACLs send higher bids than the "OFF" ones, which makes the "ON" ones always have the priority of maintaining ON state. In this case, the switching frequency can be reduced to the minimum. If $\rho=0$, the bid price is solely decided by SOA. The ACLs with the same SOA but different operation state are treated equal. It's the condition with highest switching frequency.

When the ACL is in "LockON/LockOFF" state, the ACL must keep on/off. In this case, $S_{offset}$ is assigned a large positive/negative value (here it's assigned 3/-3) to make $p_{bid}$ higher/lower than that of the ACLs in other states and ensure that the ACL has the highest priority of maintaining on/off.

### C. Virtual Market Clearing

Fig. 1(c) shows the market clearing strategy. The virtual market sorts the bids from the ACLs in descending order of the bid price and forms the demand curve. The vertical part of the supply curve equals to the target power $P_{AC}^*$ and the bid price of it is limited to ±2 to ensure that the clearing price is within the range of [-2,2]. The virtual market solves the clearing price $p^*$ at the intersection of the demand curve and the supply curve.

### D. Response to the Clearing Result

Each ACL receives and responds to the broadcast clearing price $p^*$. To ensure that the aggregated power of the ACL cluster can track the target power, each ACL should be off if its bid price is lower than $p^*$, otherwise it should be on.

Because the clearing price varies in [-2,2], which is always lower than the bid price of the "lockON" ACLs and higher than that of the "lockOFF" ones. It ensures that the "locked" ACLs can keep the current state and the lockout time constraint can be strictly satisfied.

## IV. SIMULATION RESULTS

Two cases are studied in this paper. In both cases, the performance is demonstrated with a set of 500 ACLs. The second-order ETP model[6] is adopted for modeling the ACLs. The baseline load estimation adopts the method in [6] and the estimated value is averaged to an hourly load profile.

The outdoor temperature is shown in Fig. 2.

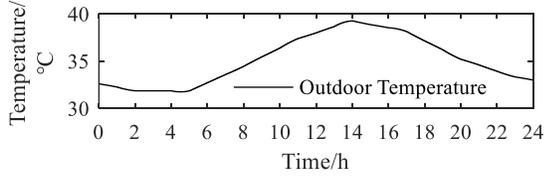

Fig. 2 Outdoor Temperature

Main parameter settings are shown in TABLE I.

TABLE I. *Main Parameter Settings*

| Area/ m² | Air Change Freq/(Times/h) | Window-Wall Ratio | SHGC | EER |
|---|---|---|---|---|
| U(88,176) | N(0.5,0.06) | N(0.15,0.01) | U(0.22,0.5) | U(3,4) |
| $R_{th}$ of Roof/ (°C.m²/W) | $R_{th}$ of Wall/ (°C.m²/W) | $R_{th}$ of Floor/ (°C.m²/W) | $R_{th}$ of Window/ (°C.m²/W) | $R_{th}$ of Door/ (°C.m²/W) |
| N(5.28,0.70) | N(2.99,0.35) | N(3.35,0.35) | N(0.38,0.03) | N(0.88,0.07) |
| $\tau$/min | $T_{desired}$/°C | $T_{high}$/°C | $T_{low}$/°C | $t_{lock}$/min |
| 30 | N(26,0.5) | U(2,3) | U(2,3) | 5 |

Note: $R_{th}$ denotes thermal resistance, $T_{high}$= $T_{max}$- $T_{set}$, $T_{low}$= $T_{set}$- $T_{min}$, U(a,b) denotes the uniform distribution and N(avg,std) denotes the normal distribution.

### A. Case 1: Mitigation of the Power Fluctuation

The uncontrollable load and the wind power are shown in Fig. 3. The control cycle is 1min.

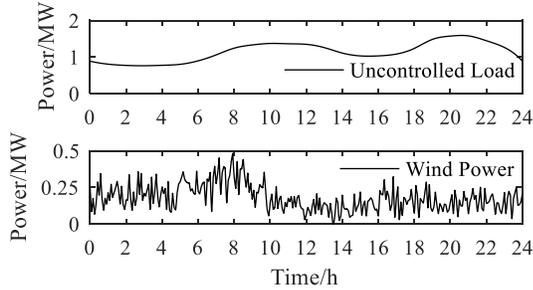

Fig. 3 Uncontrolled load and the wind power

The mitigation performance and the fluctuation rate, which adopts the definition in [6], with different $\rho$ are shown in Fig. 4, where $P_{g0}$ is the tie-line power when ACLs are not controlled. As is shown in Fig. 4, the fluctuations can be effectively mitigated with all the different $\rho$ and the response power is almost the same.

The SOA with different $\rho$ is shown in Fig. 5 with different colors and the corresponding average SOA is plotted with thicker lines. As we can see, the coverage width of SOA becomes larger when $\rho$ increases. It shows that with larger $\rho$, the ACLs are less likely to switch the operation state and the temperature deviation from the desired setpoint increases. The average SOA is almost the same with different $\rho$, indicating that $\rho$ does not affect the total energy injection caused by the external control signal. Due to the delay effect of LPF, when the tie-line power increases, the actual target power is smaller than the desired power, which makes the average SOA increase, and vice versa.

Fig. 6 illustrates the impact of $\rho$ on the daily switching cycles. To evaluate the influence of the external control, the statistics when the ACLs are uncontrolled (with the normal thermostat deadband ±1°C) is also plotted. When $\rho$ is larger than 0.5, the switching frequency is much lower than the normal condition, which indicates that if the ramping rate of the target signal is not very large, participating in such ancillary services would not lead to more mechanical wear with proper control parameters.

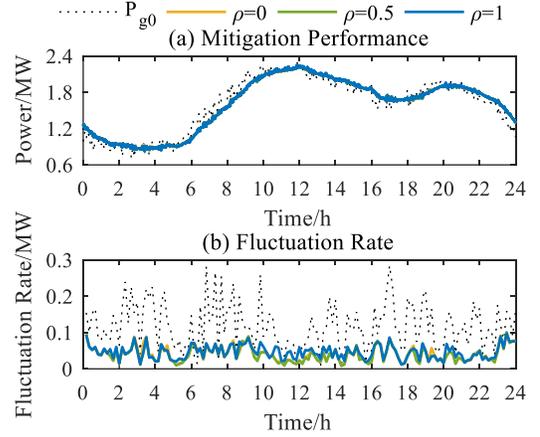

Fig. 4 The mitigation performance and fluctuation rate

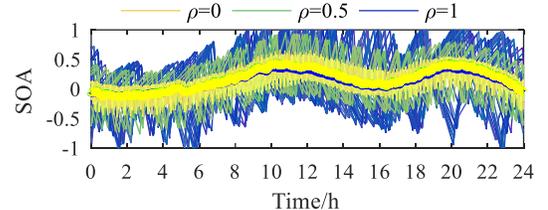

Fig. 5 The SOA of each ACL and the average SOA

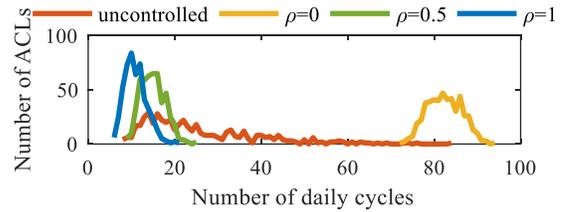

Fig. 6 The distribution of the number of daily switching cycles

In case 1, the value of $\rho$ has little impact on the control effect but large impact on the comfort level of the users and the switching frequency. Considering these two factors, it's better to set $\rho$ to the range of (0.3,0.7).

### B. Case 2: Response to Dynamic Regulation Signal

The regD signal uses the data downloaded from the PJM website[9]. The control cycle is 4s. The regulation capacity is set to 0.4MW. Due to the page limitations, only the results with two different $\rho$ are shown.

Fig. 7 shows the response performance. The grey area is the "locked" power which is the sum of the rated power of the

"lockON/OFF" ACLs, and the yellow area is the "available" power which is the sum of the rated power of the "ON/OFF" ACLs. As we can see, in this case, the response performance is affected by the value of $\rho$. The regD signal changes very fast, so when $\rho$ is too small, the operation state has much more chance to be switched than the large $\rho$ condition. It makes the "available" power too small to track the regulation signal well. When the $\rho$ is large, the "available" power can be large enough to track the target.

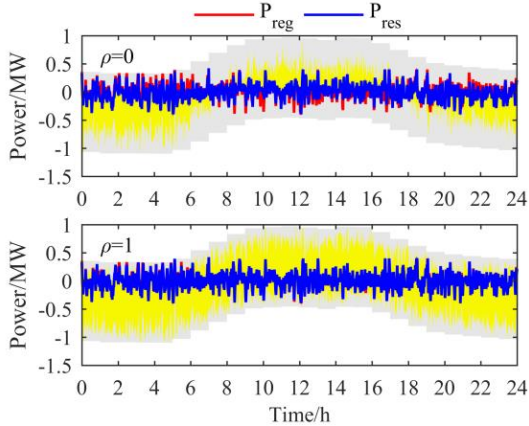

Fig. 7 The response performance

The average value of regD signal is near zero, which means the average energy injection is near zero and the average SOA can keep around zero if the baseline load is well estimated. In Fig. 8, we can see that when $\rho$ equals to 0, the average SOA deviates from the ideal value. It's because the ACLs are unable to well track the regD signal. During 0~7h and 18~24h, the ACLs consumes more power than the target, so the SOA falls down. During 8~17h, the ACLs consumes less power than the target, so the SOA rises up. In contrast, when $\rho$ equals to 1, the average SOA is kept around zero.

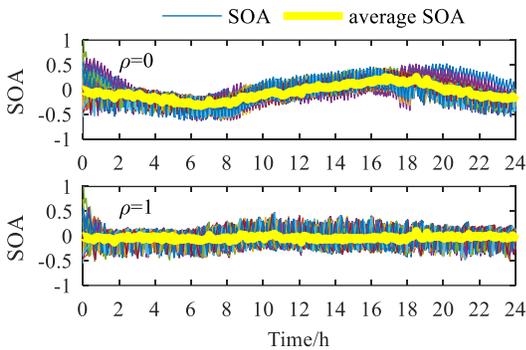

Fig. 8 The SOA of each ACL and the average SOA

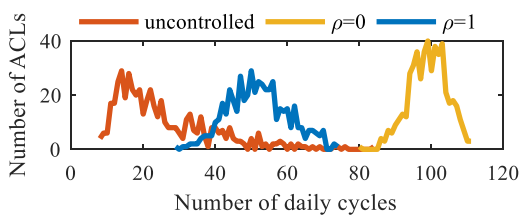

Fig. 9 The distribution of the number of switching cycles per day

In case 2, because the regD signal fluctuates too fast, even if $\rho=1$, which is the slowest condition possible, the switching frequency is still higher than the normal condition. It may exert negative impact on the devices. To enable the ACLs to participate in the fast regulation services, new control methods need to be designed to lower the switching frequency and protect the devices.

In case 2, the value of $\rho$ is recommended to set to (0.5,1) to ensure satisfactory response performance and relatively low switching frequency.

## V. CONCULSION

This paper proposes modified control strategies based on MBC method for the ACL cluster to provide ancillary services. The control framework can satisfy users' different requirements, has low communication bandwidth and can protect users' privacy. In addition, the switching frequency can be adjusted and the lockout time constraint can be strictly satisfied without increasing the computation complexity and the communication traffic.

The control parameter $\rho$ affects the response performance, switching frequency and comfort level. According to the characteristics of the target signal, $\rho$ can be properly chosen to achieve the balance between these three factors. Generally, the faster and greater the signal changes, the larger $\rho$ is required to meet the performance requirements and lower switching frequency. When the ACLs responds to relatively slowly-changing signals, smaller $\rho$ is recommended to improve the users' comfort level.

The quick response ability is the advantage of the ACLs in providing fast regulation services. However, considering the switching frequency and communication burden, innovative control strategies need to be designed. It's one of our future research direction. The reward allocation mechanism is also one research priority in later studies.